\documentclass[12pt]{article}
\usepackage[cmtip,arrow]{xy}
\usepackage{pb-diagram, pb-xy}
\usepackage{amssymb,epsfig,amsfonts}
\newtheorem{theorem}{Theorem}
\newtheorem{coroll}{Corollary}
\newtheorem{lemma}{Lemma}
\newtheorem{prop}{Proposition}
\newcommand{\ra}{\longrightarrow}
\newcommand{\da}{\downarrow}
\newcommand{\cO}{\hat {\cal O}}
\newcommand{\cF}{\hat {\cal F}}
\newcommand{\pr}{\prod\limits}
\def\R{{\mathbb R}}
\def\Z{{\mathbb Z}}

\def\C{{\mathbb C}}

\def\A{{\mathbb A}}
\author{}
\date{}
\begin{document}

\vspace*{2cm}
\begin{center}
{\bf INTEGRABLE SYSTEMS AND LOCAL FIELDS} \\
\par\bigskip
{\sc A. N. Parshin} \\
\end{center}
\par\bigskip
\begin{flushright}
{\em To the memory of \\
Alexey Ivanovich Kostrikin}
\end{flushright}

\par\bigskip
In 70's there was discovered a construction how to attach to some
algebraic-geometric data an infinite-dimensional subspace in the space
$k((z))$ of the Laurent power series. The construction was
successfully used in the theory of integrable systems, particularly, for
the KP and KdV equations \cite{Kr, SW}.
There were also found some applications to the moduli of algebraic curves
\cite{ADKP, BS}.
Now it is known as the Krichever correspondence or the Krichever map
\cite{ADKP, M, AMP, Q, BF}. The original work by I. M. Krichever has also
included commutative rings of differential operators as a third part of
the correspondence.

The map we want to
study here was first described in an explicit way by G. Segal and G. Wilson
\cite{SW}. They have used an analytical version of the infinite dimensional
Grassmanian introduced by M. Sato \cite{SS, PS}. In the sequel we consider a
purely algebraic approach as developed in \cite{M}.

Let us just note that
the core of the construction is an embedding of the affine coordinate ring
on an algebraic  curve into the field $k((z))$ corresponding to the power
decompositions in a point at infinity
(the details see below in section 2). In number theory this corresponds to
an embedding of the ring of algebraic integers to the fields $\C$ or $\R$.
The latter one is well known starting from the XIX-th century. The  idea
introduced by Krichever was to insert the local parameter $z$. This
trick looking so simple enormously extends the area of the correspondence.
It allows to consider all algebraic curves simultaneously.

But there still remained a hard restriction by the case of curves, so by
dimension 1. Recently, it was pointed out by the author \cite{P3,P4} that there
are some connections between the theory of the KP-equations and the theory of
$n$-dimensional local fields \cite{P1}, \cite{FP}. From this point of view
it becomes clear that the Krichever construction should have a generalization
to the case of higher dimensions. This generalization is suggested in the
paper  for the case of algebraic surfaces (see theorem 4 in section 4). A
further generalization  to the case of arbitrary dimension was recenly
proposed by D. V. Osipov \cite{O}.

We start with  description of a connection between the KP hierarchy
in the Lax form and the vector fields on infinite Grassmanian manifolds
(section 1). These results are known but we prove them here in more
transparent and simple way (we have used basically \cite{M} and \cite{LM}).
In appendix 1 we remind how to get the standard
KP and KdV equations from the Lax operator form. In appendix 2 we ouline a construction
of the semi-infinite monomes for the field $k((z))$ which is an important part
of the theory of Sato Grassmanian.

An important feature of all these considerations is their purely
algebraic character. Everything can be done over an arbitrary field
$k$ of characteristic zero\footnote{Further developement of different aspects of
the Krichever correspondence
for dimension two see in \cite{ZhO, KOZh1,  KOZh2}}.

Let us also note that the construction of the restricted adelic complex
in section 3 is of an independent interest, also in arithmetics. It has
already appeared in a description of vector bundles on algebraic surfaces
\cite{P2}.
\par\medskip
This work was partly done during my visit to the
Institut f\"ur Algebra
und Zahlentheorie der Technische Universit\"at Braunschweig in 1999.
I am very much grateful to Professor Hans Opolka for the hospitality\footnote{The text was published in
{\em Communications in Algebra}, 29(2001), No.9, 4157-4181. This version includes a corrected
proof of the proposition 2. I'm grateful to Alexander Zheglov for the correction of a mistake in the original proof.
Also,  we include some additional remarks on the deduction of concrete equations from
the Lax hierarchy and appendix 2.}.

\section{Sato correspondence for dimension 1}

This correspondence connects two seemingly distant objects: infinite
Grassmanian manifolds and rings of pseudo-differential operators.

Let $K = k((z))$ be the field of Laurent power series with
 filtration  $K(n) = z^nk[[z]]$. If $V \cong k((z))$ is a vector
 space of dimension 1 over $K$ then we can choose a filtration $V(n)$
 such that $K(n)V(m) \subset V(n+m)$. Let $V_1 : = V(0)$.

Denote by $\mbox{Gr}(V)$
the set of the subspaces $W$ in $V$ such that the complex
$$
W \oplus V_1 \rightarrow V\eqno (1)
$$
is a Fredholm one.
It is a
(infinite-dimensional) projective variety with connected components
marked by the Euler characteristic of the complex (1)~\cite{AMP,K}(see also, appendix 2).

     Let us now introduce the ring
     $E = k[[x]]((\partial^{-1}))$ of formal pseudo-differential operators with
     coefficients from the ring $k[[x]]$ of regular formal power series as
     the left $k[[x]]$-module of all   expressions
     $$ L = \sum_{i > -\infty}^{n} a_i \partial^i,~a_i \in k[[x]]. $$
     Then a multiplication can be defined according to the Leibnitz rule:
     $$ (\sum_i a_i\partial^i)(\sum_j b_j\partial^j)      \sum_{i,j;k \ge 0} {i\choose k}a_i d^k(b_j)\partial^{i + j - k}.$$
     Here we put
     $${i\choose k} = \frac{i(i-1)\dots (i-k+1)}{k(k-1)\dots 1},
     ~\mbox{if}~k > 0,~{i\choose 0} = 1$$
     and $d$ is the derivation by $x$.

     Particularly, for  $f \in k[[x]]$:
     $$[\partial,f] = \partial f - f\partial = d(f),$$
      (Heisenberg commutative relation),
     $$[\partial^{-1},f] = \partial^{-1} f - f \partial^{-1} -d(f)\partial^{-2} + d^2(f)\partial^{-3}
     - \dots .$$
     It can be checked that $k[[x]]((\partial^{-1}))$ will be
an associative
     ring (the details see in \cite{P4}).

There is a decomposition
$$ E = E_{+} + E_{-},$$
where $E_{-} = \{L \in E: L = \sum_{n < 0} a_n \partial^n \}$~and
~$E_{+}$~consists of the operators containing only $ \ge 0$~powers~of
~$\partial$.
The elements from $E_{+} =: D$ are the differential operators and the elements
from $E_{-}$ are the Volterra operators.

From the commutation relations we see that  $E = k((\partial^{-1})) \oplus
Ex$ and thus the map $E \rightarrow E/Ex = V$ (we identify the image
of $\partial^{-1}$
with $z$) defines a linear action of the ring $E$ on  $V$ and also on
$\mbox{Gr}(V)$. The subspace $V_1$ is transformed by the action of operators
from $E$ into a subspace $V^{'}$ commesurable with $V_1$ (it means that
the quotient $V^{'}+V_1/V^{'} \cap V_1$ is of finite dimension). Thus the
Fredholm condition from the definition of Grassmanian manifold will be
preserved.

We will call the map $E \rightarrow V$ by the Sato map.

\begin{prop}({\sc Lemma on a stabilizer}). Let $P \in E$ and
$W_0 = k[z^{-1}]  \in Gr(V)$ if $V = k((z))$.

Then  $PW_0 \subset W_0$ if and only if $P \in D$.
\end{prop}
{\sc Proof}. Since $W_0$ is equal to the image of $E_{+}$ under the Sato
map (by lemma 1 below) we can replace the first condition by the following one
$$
P\cdot E_{+}  \subset E_{+} + Ex\eqno (2)
$$
and work in the ring $E$. Since $E_{+}$ is a ring we have $P\cdot E_{+}
\subset E_{+}$ for $P \in E_{+} = D$.

Now assume (1) and decompose $P$ as $P = P_{+} + P_{-}$. We will also use
the notation $L \sim M$ for $L, M \in E$ such that $L - M \in Ex$.

\begin{lemma}. In the ring $E$, $x^n\partial^n \sim c_n$ where
$c_n = (-1)^n n!$.
\end{lemma}
{\sc Proof}. By commutation relations $x^n\partial^n = \partial^n x^n +
(\mbox{some coefficient})\partial^{n-1}x^{n-1} + \dots + c_n$ and applying this
operator to $x^{-1}$ we get the value of $c_n$.
\par\smallskip
Returning to our proposition we see that (2) implies
$$
P_{-}E_{+} \in E_{+} + E_{-}x.\eqno (3)
$$
To prove that $P_{-} = 0$ it is enough to show $P_{-} \in E_{-}x^n$ for
all $n \ge 1$.

First, $1 \in E_{+}$ and consequently $P_{-} \in E_{+} + E_{-}x E_{+} \oplus E_{-}x$. We see that $P_{-} \in E_{-}x$.

Proceeding by induction we assume $P_{-} = Q_{-}x^n$ with $Q_{-} \in E_{-}$.
We have
$$
Q_{-}x^n\partial^n \in E_{+} + E_{-}x,~\mbox{by}~(3)
$$
$$
Q_{-}x^n\partial^n \in Q_{-}Ex + Q_{-} \subset Ex + Q_{-},~\mbox{by lemma}.
$$
Taking all together we get $Q_{-} \in  E_{+} + Ex \subset  E_{+} \oplus E_{-}x$.
It means $Q_{-} \in  E_{-}x$ and $P_{-} \in E_{-}x^{n+1}$.

\begin{prop}({\sc Transitivity theorem}). Let $W \in Gr(V)$ and $W \oplus V_{1}
=V$. Then there exists an unique operator $S \in 1+ E_{-}$ such that $W = S^{-1}W_0$.
\end{prop}
{\sc Proof}. If $W$ satisfies the conditions of the theorem then $W$ is a
union of subspaces $W \cap V_{-n}$ and one can choose basis $w_n$ in $W$ such
that $w_n = z^{-n}  + v_n,~v_n \in V_1$. We want to construct an operator
$S = 1 + P \in 1 + E_{-}$ such that $Sz^{-n} = w_n+c_0w_0+\ldots +c_{n-1}w_{n-1}$ for some $c_0,\ldots c_{n-1}\in k$, or $Pz^{-n} = v_n+c_0w_0+\ldots +c_{n-1}w_{n-1}$ for
all $n \ge 0$.

The Sato map $E \rightarrow V$ transforms $\partial^n$ into $z^{-n}$ and we
can find $Q_n \in E_{-}$ such that $Q_n$ goes to $v_n$ by the map. Thus in
order  to construct  our $S$ from the conclusion of the theorem we have to
find an operator $P \in E_{-}$ such that
$$
P\partial^n = Q_n +c_0(1+Q_0)+\ldots +c_{n-1}(\partial^{n-1}+Q_{n-1})+ Ex,\eqno (4)
$$
where $Q_n$ is a given sequence of operators from $E_{-}$. From now on we
work in the ring $E$. Put $P_0=Q_0$.

{\sc CLAIM}. For $n \ge 1$ there exist operators $P_n \in E_{-}$ such that
\begin{quotation}
i) $ P_n\partial^n+ P_0\partial^n + \dots + P_{n-1}\partial^n  \in Q_n +c_0(1+Q_0)+\ldots +c_{n-1}(\partial^{n-1}+Q_{n-1})+ Ex,~\mbox{for some }c_0,\ldots , c_{n-1}\in k $

ii) $P_n, P_n\partial, \dots, P_n\partial^{n-1} \in Ex$

iii) $P_n \in Ex^n$

\end{quotation}

First, we show that the claim implies the existence of the $P$ with the
property (4).
We put $P = P_0 + P_1 + \dots$. Then
$$
\begin{array}{ccc}
P_0\partial^n + P_1\partial^n  + \dots + P_n\partial^n + P_{n+1}\partial^n +
\dots & \sim & \mbox{(by ii)} \\
P_0\partial^n + P_1\partial^n  + \dots +
P_n\partial^n & \sim & \mbox{(by i)} \\
Q_n +c_0(1+Q_0)+\ldots +c_{n-1}(\partial^{n-1}+Q_{n-1}).
\end{array}
$$
Here we use the notation $\sim$ from the lemma 1. The property iii) implies
the convergence of the series for $P$ by the following
\begin{lemma}. In the ring $E$, the series $\sum_{n \ge 0} P_n, P_n \in
E_{-}$ will converge if $P_n \in E_{-}x^n$.
\end{lemma}
Now we prove the claim and then the lemma.

We can define $P_n$ by induction. Obviously, by induction we have
$$
P_0\partial^n + \dots + P_{n-1}\partial^n\in c_0(1+Q_0)+\ldots +c_{n-1}(\partial^{n-1}+Q_{n-1})+ Ex +E_-
$$
for some $c_0,\ldots , c_{n-1}\in k $. So, the image of the operator
$P_0\partial^n + \dots + P_{n-1}\partial^n-c_0(1+Q_0)+\ldots +c_{n-1}(\partial^{n-1}+Q_{n-1})$ under the Sato map lies in $V_1$ and we can find an operator $Q'_n\in E_-$ such that $Q'_n$ goes to this image by the map.

Then we can take $P_n = (Q_n-Q'_n)x^n\in E_-$. It gives iii) and
$$
P_n\partial^n = (Q_n-Q'_n)x^n\partial^n \sim Q_n-Q'_n \sim
$$
$$
Q_n +c_0(1+Q_0)+\ldots +c_{n-1}(\partial^{n-1}+Q_{n-1})-P_0\partial^n - \dots - P_{n-1}\partial^n
$$
by lemma 1 and
$$
P_n\partial^k = Q_n^{'}x^n\partial^k = Q_n^{'}x^{n-k}x^k\partial^k \sim
Q_n^{'}x^{n-k} \sim 0,~k = 0,1,\dots,n-1
$$
again by lemma 1 and we are done.

The uniqueness of the operator $S \in 1 + E_{-}$ such that $W = S^{-1}W_0$
follows
from proposition 1. Indeed, let $SW_0 = W_0$. Then  $S$ must belong to
$E_{+}$ and thus $S = 1$.

{\sc Proof} of Lemma 2. Let $P_n = P_n^{'}x^n$ where
$$
P_n^{'} = \sum_{m \ge 0} a_m^{(n)}\partial^{-m}~\mbox{and}~
P_n  = \sum_{m \ge 0} b_m^{(n)}\partial^{-m}.
$$
We see $b_m^{(n)} = a_m^{(n)}x^n \pm a_{m-1}^{(n)}d(x^n) \pm a_{m-2}^{(n)}
d^2(x^n) \pm \dots$
and for every $m$ the series $\sum_{n \ge 0} b_m^{(n)}$ will converge
in $k[[x]]$.
\par\bigskip
If $W \in \mbox{Gr}(V)$ then  we have
$T_W = \mbox{Hom}(W, V/W)$ for the tangent space in the point $W$
and there is a natural map
$\mbox{Hom}(V, V) \rightarrow T_W$. For $n \in \Z$  we can define
a vector field $T_n$ on $\mbox{Gr}(V)$. It is equal to the image of
the multiplication operator by $z^{-n}$ in the space $V$.

The KP hierarchy is a dynamical system defined on an affine space
$E':=  \partial + E_{-}$. The tangent space to any point $L \in E'$
is canonically equal to $E_{-}$.
\par\smallskip
{\sc Definition 1}. The $n$-th vector field of the KP hierarchy on $E'$
is defined as $KP_n = [(L^n)_{+}, L]$.
\par\smallskip
Since $(L^n)_{+} = L^n - (L^n)_{-}$ the field $KP_n$ belongs to $E_{-}$.

The set $G = 1 + E_{-}$ is a group carrying the vector fields
$-(S\partial^n S^{-1})_{-}S$. The tangent space to any point from $G$ is
again $E_{-}$.

At last we denote by $Gr_{+}$ the big cell from the Grassmanian manifiold,
$$
Gr_{+}(V) = \{W \in Gr(V): W \oplus V(1) = V \}.\eqno (5)
$$

All these spaces, $E', G, Gr_{+}$ have the distinguished points: $\partial,
1, W_0 = k[z^{-1}]$(if $V = k((z))$).

\begin{theorem}({\sc Sato correspondence}). The maps
$$
E'\stackrel{\varphi}{\leftarrow}G\stackrel{\psi}{\rightarrow}\mbox{Gr}_{+}(V),
$$
where $\varphi(S) = S\partial S^{-1}, \psi(S) = S^{-1}(W_0)$
have the following properties
\begin{quotation}
i) for any $S \in G$, $L = \varphi(S) \in E'$ and $W = \psi(S) \in Gr_{+}$
the diagram
$$
\begin{array}{ccccc}         `
T_L & \stackrel{d\varphi_S}{\leftarrow} & T_S & \stackrel{d\psi_S}{\rightarrow}
& T_W \\
&  & & & \parallel \\
\parallel &  & \parallel && Hom(W, V/W) \\
&  & & & \uparrow \\
E_{-} & \stackrel{\varphi '}{\leftarrow} & E_{-}
& \stackrel{\psi '}{\rightarrow}& Hom(V, V)
\end{array}
$$
commutes. Here $d\varphi_S$ and $d\psi_S$ are jacobian maps of the maps
$\varphi,\psi$ on the tangent space $T_S$, and
$$
\varphi '(A) = [AS^{-1}, L],
$$
$$
\psi '(A) = -S^{-1}A~\mbox{acting on}~V~\mbox{by the Sato map}
$$
with $A \in E_{-}$.

ii) for any $S \in G$, $L = \varphi(S) \in E'$ and $W = \psi(S) \in Gr_{+}$
and any $n \ge 1$
$$
d\varphi_S(-(S\partial^n S^{-1})_{-}S) =  KP_n,
$$
$$
d\psi_S(-(S\partial^n S^{-1})_{-}S)  = T_{n}
$$
\end{quotation}
\end{theorem}

{\sc Proof}. First we consider the left hand side of the diagram. If
$A \in E_{-} = T_S$ and $R = S + A$ is an infinitely small deformation of
$S$ then  up to the higher powers of $A$
$$
R = S(1 + S^{-1}A) =  (1 + AS^{-1})S,
$$
$$
R^{-1} =  (1 - S^{-1}A)S^{-1} = S^{-1}(1 - AS^{-1}),
$$
$$
R\partial R^{-1} = (1 + AS^{-1})S\partial S^{-1}(1 - AS^{-1}) $$
$$
(S\partial S^{-1}  + A\partial S^{-1})(1 - AS^{-1}) $$
$$
S\partial S^{-1} + A\partial S^{-1} - S\partial S^{-1}AS^{-1} $$
$$
L + [AS^{-1}, L] = L + d\varphi_S(A).
$$
Using this result we can check up the statement on the vector fields:
$$
d\varphi_S(-(S\partial^n S^{-1})_{-}S) = [-(S\partial^n S^{-1})_{-}SS^{-1},
L] = [(L^n)_{+}, L].
$$
\par\smallskip
It remains to consider the right hand side of the diagram. Let $W = S^{-1}
W_0$ as in the diagram and let $R = S + A$ be as above. Then $W': R^{-1}W_0$ and we get
$$
W' = R^{-1}SS^{-1}W_0 = R^{-1}SW = (1 - S^{-1}A)S^{-1}SW = (1 - S^{-1}A)W.
$$
Since $W, W'$ are two spaces from the big cell
$$
W \oplus V(1) = V,
$$
$$
W' \oplus V(1) = V
$$
and the space $W'$ defines a linear map
$$
W \rightarrow V \rightarrow V/W' \stackrel{\sim}{\leftarrow} V(1)
\stackrel{\sim}{\rightarrow} V/W
$$
which is an element of the tangent space $T_W$ corresponding to the
deformation $W'$ of $W$ (see, for example, \cite{K}). It is easy to see that for our space $W'$
the linear map will coincide with the action of operator $-S^{-1}A$
through the Sato map.
\par\smallskip
The last step of the proof is to check that $d\psi_S$ takes
$-(S\partial^n S^{-1})_{-}S$ into $T_{n}$. But we have
$$
-S^{-1}(-(S\partial^n S^{-1})_{-}S) = S^{-1}L^nS - S^{-1}(L^n)_{+}S \partial^n - S^{-1}(L^n)_{+}S
$$
and we have to show that the second term is trivial in $T_W$.

This can be seen from the commutative diagram
$$
\begin{array}{cccccccc}
& W_0 & \rightarrow & V & \rightarrow & V & \rightarrow  & V/W_0 \\
S^{-1} & \downarrow &  S^{-1} & \downarrow &
S^{-1} & \downarrow &  S^{-1} & \downarrow \\
& W & \rightarrow & V & \rightarrow & V & \rightarrow  & V/W
\end{array}
$$
The bottom map from $V$ to $V$ is equal to the Sato image of
$S^{-1}(L^n)_{+}S$ and
the corresponding top horisontal map is the Sato action of the operator
$(L^n)_{+}$.
By the proposition 1 the last one is trivial in $T_{W_0} = Hom(W_0, V/W_0)$.

The theorem is proved.

{\sc Remark 1}. The maps $\varphi$ and $\psi$ can be deduced from the corresponding
actions of the group $G$ on the manifolds $E'$ and $Gr(V)$. One has to
consider the actions on the orbits going through $\partial$ and $W_0$,
respectively. The first orbit is one of the co-adjoint orbits for the
(infinite dimensional) Lie group $G$.

\begin{coroll}. The Sato correspondence induces
the diagram of bijections
$$
E'\stackrel{\varphi}{\leftarrow}G/G_0
\stackrel{\psi}{\rightarrow}\mbox{Gr}_{+}(V)/k[[z]]^* ,
$$
where $G_0 = G \cap k((\partial^{-1}))$ and the action of $k[[z]]^*$ on
$Gr(V)$ is defined by the module structure on $V$ over $K$.
\end{coroll}
{\sc Proof} for the map $\psi$ easily follows from the definitions, theorem 1
and proposition 2. To check up the bijectivity of the map $\varphi$ one has
to apply theorem 1 from \cite{P4}.

\section{Krichever correspondence for dimension 1}
We first discuss the adelic complexes for the case of dimension 1.
Concerning a definition of the adelic notions we refer to \cite{FP},\cite{H}.
We also note that the sign $\prod$ denotes the adelic product.

Let $C$ be an projective algebraic curve over a field $k$, $P$ be
a smooth point  and $\eta$ a general point on $C$.
For every point (in Grothendieck's sense) $\alpha \in C$ we have a field
$K_{\alpha}$. $K_{\alpha}$ is a quotient ring of the completed local
ring ${\cO}_{\alpha}$ of the point $\alpha$. It is a 1-dimensional local
field.

Let ${\cal F}$ be a
torsion free coherent sheaf on $C$. We denote by $\cF_{\alpha} = {\cal F} \otimes
\cO_{\alpha}$ the completed fiber at the point $\alpha \in C$. The fiber
${\cal F}_{\eta}$ at  a general point is also the space of all rational
sections of the sheaf ${\cal F}$.
\begin{prop}. The following complexes are quasi-isomorphic:
\begin{quotation}
i) adelic complex
$$ {\cal F}_{\eta} \oplus \pr_{x \in C} \cF_x \ra \prod_{x \in C}
(\cF_x \otimes_{\cO_x} K_x)
$$
ii) the complex
$$ W \oplus \cF_P \ra \cF_P \otimes_{\cO_P} K_P
$$
where  $W = \Gamma(C - P, {\cal F}) \subset \cF_{\eta}$.
\end{quotation}
\end{prop}
{\sc Proof} will be done in two steps. First, the adelic complex
contains a trivial exact subcomplex
$$  \pr_{x \in U} \cF_x \ra \pr_{x \in U} \cF_x,
$$
where $U = C - P$. The quotient-complex is equal to
$$
{\cal F}_{\eta} \oplus \cF_P \ra \prod_{x \in U}
(\cF_x \otimes K_x)/\cF_x  \oplus \cF_P.
$$
It has a surjective homomorphism to the exact complex
$$
{\cal F}_{\eta}/W \ra \prod_{x \in U}
(\cF_x \otimes K_x)/\cF_x.
$$
The exactness of the complex is the strong approximation theorem
for the curve $C$ \cite{B}[ch.II, \S 3, corollary of prop. 9; ch. VII,
\S, prop. 2]. The kernel of this
surjection will be the second complex from proposition.
\par\smallskip

 Let us now  explain the Krichever correspondence for dimension 1.

 {\sc Definition 2}.
 $$
 \begin{array}{lll}
 {\cal M}_1 & := & \{ C, P, z, {\cal F}, e_P \}\\
 C & & \mbox{projective irreducible curve}~/k \\
 P \in C && \mbox{a smooth point} \\
 z && \mbox{formal local parameter at}~P \\
 {\cal F} && \mbox{torsion free rank}~$r$~\mbox{sheaf on}~C \\
 e_P && \mbox{a trivialization of}~{\cal F}~\mbox{at}~P
 \end{array}
 $$
 Independently, we have the field $K = k((z))$ of Laurent power series with
 filtration  $K(n) = z^nk[[z]]$. Let $K_1 : = K(0)$. If $V = k((z))^{\oplus r}$
 then $V(n) = K(n)^{\oplus r}$ and $V_1 : = V(0)$.

\begin{theorem}\cite{M}. There exists a canonical map
$$
\Phi_1 : {\cal M}_1 \ra \{\mbox{vector subspaces}~A \subset K, W \subset V \}
$$
such that
\begin{quotation}
i) the cohomology of complexes
 $$
A \oplus K_1 \ra K, ~W \oplus V_1 \ra V
 $$
are isomorphic to $H^{\cdot}(C, {\cal O}_C)$ and $H^{\cdot}(C, {\cal F})$,
respectively

ii) if $(A, W) \in \mbox{Im}~\Phi_1$ then $A \cdot A \subset A,
A \cdot W \subset W$,

iii) if $m,m^{\prime} \in {\cal M}_1$ and $\Phi_1(m) \Phi_1(m^{\prime})$ then $m$ is isomorphic to $m^{\prime}$
\end{quotation}
\end{theorem}

{\sc Proof}. If~$m = (C, P, z, {\cal F}, e_P) \in {\cal M}_1$ then we put
$$
A : = \Gamma(C - P, {\cal O}_C),
$$
$$
W: = \Gamma(C - P, {\cal F}).
$$
Also we have by the choice of $z$ and $e_P$
$$
\cO_P = k[[z]],~K_P = k((z)),
$$
$$
{\cal F}_P  = {\cal O}_Pe_P = {\cal O}_P^{\oplus r},~\cF_P = \cO_P^{\oplus r}.
$$
This defines the point $\Phi_1(m) \in {\cal M}_1$. Indeed, for the subspace
$W$ we have the following canonical identifications
$$
\Gamma(C - P, {\cal F}) \subset {\cal F}_{\eta} \otimes_{{\cal O}_P} K_P
= \cF_P \otimes K_P = \cO_P^{\oplus r} \otimes K_P = k((z))^{\oplus r}.
$$
The same works for the subspace $A$.

The property ii) is
obvious, the property i) follows from the proposition 3. To get iii) let
us start with a point $\Phi_1(m) = (A, W)$. The standard valuation on $K$
gives us  increasing filtrations $A(n) = A \cap K(n)$ and
$W(n) = W \cap V(n)$ on the spaces $A$ and $W$. Then we have
$$
\begin{array}{ccl}
C - P & = & \mbox{Spec}(A),  \\
C  & = &  \mbox{Proj}(\oplus_n A(n)), \\
{\cal F} & = & \mbox{Proj}(\oplus_n W(n)),
\end{array}
$$
by lemma 9. Thus we can reconstruct the quintiple $m$ from the point
$\Phi_1(m)$.

{\sc Remark 2}. It is possible to replace the ground field $k$ in the
Krichever construction by an arbitrary scheme $S$, see \cite{Q}.

Using the Krichever correspondence $\Phi_1$ one can construct the integral
varieties in $X = \mbox{Gr}(V)/k[[z]]^*$ for the vector fields $T_n$ (see
section 1).  Let us fix $C, P, z$. Then the
image $\Phi_1(C, P, z, {\cal F}, e_P)$
does not depend on $e_P$ in $X$ and will run through the generalized jacobian
$\mbox{Jac}(C)$ of the curve  $C$ when we vary the invertible sheaf
${\cal F}$.

To show this fact we consider the commutative diagram
$$
\begin{array}{ccccccccc}
&&&& k[[z]] & = & k[[z]] && \\
&&&& \downarrow && \downarrow &&\\
0 & \rightarrow & A_W & \rightarrow & K & \stackrel{\alpha}{\rightarrow} &
Hom(W, V/W) & = & T_{Gr, W}\\
&& \parallel && \downarrow && \downarrow && \\
0 & \rightarrow & A_W & \rightarrow & K/k[[z]] & \rightarrow & T_{X, W} &&\\
&&&& \downarrow && \downarrow && \\
&&&&    0  && 0   &&
\end{array}
$$
Here  $\alpha$ is the action of $K$ on $V$ by multiplications, $W \in Gr(V)$
is the second space from the image $\Phi_1(C, P, z, {\cal F}, e_P)$
and    $A_W = \{ f \in K: fW \subset W \}$.

The diagram explains what happens with tangent spaces when we go from
the Grassmanian manifold to it's quotient $X$ by the group $k[[z]]^*$.
The space $k[[z]]$ is the Lie algebra of the last group.

For the case of invertible sheaf  ${\cal F}$ we know that $A_W = A$ (see
\cite{SW}[n 6]). Thus we get the exact sequence
$$
A + k[[z]] \rightarrow k((z)) \rightarrow T_{X,W},
$$
where $k((z))/A + k[[z]] = H^1(C, {\cal O}_C)$  is the tangent space
to the generalized jacobian of the curve $C$. From the sequence we conclude
that all vector fields $T_n$ belong to the image of $H^1(C, {\cal O}_C)$
in $T_{X,W}$ and consequently they are tangent to the image of generalized
jacobian $\mbox{Jac}(C)$ of the curve $C$ under the map $\Phi_1$ (see
details in \cite{M}).

This result works for the component of $Gr(V)~(V = k((z)))$ containing
the subspace $W$ such that $W \oplus k[[z]] = V$. The image of the Sato
correspondence belongs to another component containing the space $W_0$ and
the big cell $Gr_+(V)$ (see (5), section 1).

The multiplication by $z$ transforms one component onto another preserving
the vector fields $T_n$. Thus we get the integral varieties for the KP
flow which are abelian varieties for smooth curves $C$. If $k = \C$ they
are topological toruses with the movement along the straight lines.

This explain why the KP system can be considered to some extent as an
integrable one. But we must have in mind that there is a lot of points
$W$ in $Gr(V)$ where the dynamical behavior is quite different. Take,
for example, the points with $A_W = k$.

If $C$ is a projective line with
a double point then $\mbox{Jac}(C) = k^*$  and we get a  1-soliton solution
of the KP equation.

We see that local fields enters into this picture in   essential way.
The ring $E$ is a subring of 2-dimensional local skew-field
$P = k((x))((\partial^{-1}))$ (see \cite{P4}). The fields $k((z))$ and
$K_P$ are 1-dimensional local fields. The complex (1)(section 1) from
definition of the Grassmanian manifold intimately connected with adelic
complex on an algebraic curve (property i) from Theorem 2).

We suggest that analogous constructions should exist for
higher dimensions as well
(concerning the rings of PDO see [11]). Here we consider a generalization
of the Krichever map $\Phi_1$ to the case of algebraic surfaces.

\section{Adelic complexes in dimension 2}

Let $X$ be a projective irreducible algebraic surface over a field
     $k$, $C \subset X$ be an irreducible projective curve, and
     $P \in C$ be a smooth point on both $C$ and $X$. Let ${\cal F}$
be a torsioh free coherent sheaf on $X$.

We remind the definition of the standard rings attached to a pair
$x \in D \subset X$ on the surface $X$. Here the $D$ is an irreducible
divisor. It correspond to an ideal $\wp \subset {\cal O}_x$ in the local
ring of the point $x$. First we apply localization by $\wp$ to ${\cO}_x$
and then take a completion by the ideal $\wp$. We get a  ring ${\cal O}_{x, D}$
which is a complete discrete valuation ring if the point $x$ is smooth
on both $D$ and $X$. The local field $K_{x, D}$ is a quotient ring of the
${\cal O}_{x, D}$. It is really a field in the smooth case. In this case
we have  ${\cal O}_{x, D} = k((u))[[t]],~K_{x, D} = k((u))((t))$ if
${\cO}_x = k[[u, t]]$ and $\wp = (t)$. The field $K_{x, D}$ is an example
of 2-dimensional local field.

There are some rings attached to the point $x$ and the divisor $D$.
Denote by $K$ the field of rational functions on the $X$. Let $K_x K{\cO}_x \subset $ a quotient ring of the local ring $\cO_x$.
If $\cO_D$ is a local ring of the divisor $D$ then let $K_D$ be it's
quotient ring. In the smooth case $\cO_D = k(D)[[t]],~K_D = k(D)((t))$.

{\sc Definition 3}. Let $x \in C$. We let
     $$ B_x({\cal F}) = \bigcap\limits_{D \ne C}((\cF_x \otimes K_x) \cap
(\cF_x \otimes \cO_{x,D})),$$
     where the intersection is done inside the group $\cF_x \otimes K_{x}$,
     $$ B_C({\cal F}) = (\cF_C \otimes K_C) \cap
     (\bigcap\limits_{x \ne P} B_x),$$
     where the intersection is done inside $\cF_x \otimes K_{x,C}$,
     $$ A_C({\cal F}) = B_C({\cal F}) \cap \cF_C,$$
     $$ A({\cal F}) = \cF_{\eta} \cap (\bigcap\limits_{x \in X-C} \cF_x).$$

We will freely use the following shortcuts:
$$
\begin{array}{lll}
K\cF_x & = &  \cF_x \otimes_{\cO_x} K_x, \\
K\cF_D & = &  \cF_D \otimes_{\cO_D} K_D, \\
\cF_{x,D} & = &  \cF_x \otimes_{\cO_x} {\cal O}_{x,D}, \\
K\cF_{x,D} & = &  \cF_x \otimes_{\cO_x} K_{x,D}.
\end{array}
$$
Next, we need two lemmas connecting the adelic complexes on $X$ and $C$.
They are the versions of the relative exact sequences, see \cite{P1},
\cite{FP}.
The curve $C$ defines the following ideals:
$$
K_{x,C} \supset \cO_{x,C} \dots \supset \wp^n_{x,C} \supset \dots,
$$
$$
K_C \supset \cO_C \dots \supset \wp^n_C \supset \wp^{n+1}_C \supset \dots,
$$
$$
K_x \supset \cO_x \dots \supset \wp^n_x \dots ,
$$
and $\wp_x = \cO_x \cap \wp_{x,C}$.

\begin{lemma}. We assume that the curve $C$ is a locally complete intersection.
Let $N_{X/C}$ be the normal sheaf for the curve $C$  in $X$.
For all $n \in \Z$ the maps
$$
\prod_{x \in C} \wp^n_{x,C}\cF_{x,C}/\wp^{n+1}_{x,C}\cF_{x,C}   \ra
\A_{C,01}({\cal F} \otimes \check{N}_{X/C}^{\otimes n}),
$$
$$
\prod_{x,C} \wp^n_x\cF_{x}/\wp^{n+1}_x\cF_{x} \ra \A_{C,1}({\cal F}
\otimes \check{N}_{X/C}^{\otimes n}),
$$
$$
\wp^n_C\cF_{C}/\wp^{n+1}_C\cF_{C} \ra \A_{C,0}({\cal F}
\otimes \check{N}_{X/C}^{\otimes n}),
$$
are bijective.
\end{lemma}
In general, we have an exact sequence
$$
0 \ra  {\cal J}^{n+1}  \ra {\cal J}^{n} \ra {\cal J}^{n}\vert_C \ra 0
$$
where ${\cal J} \subset  {\cal O}_X$ is an ideal defining the curve $C$.
In our case ${\cal J} = {\cal O}_X(-C)$ and $N_{X/C} = {\cal O}_X(C)\vert_C$.
Thus the isomorphisms from the lemma are coming from the
exact relative sequence
$$
0 \ra \A_X({\cal F}(-(n+1)C)) \ra \A_X({\cal F}(-nC)) \ra
\A_C({\cal F}(-nC)\vert_C) \ra 0.
$$

\begin{lemma}. Let $P \in C$. For all $n \in {\bf Z}$  the complex
$$
\wp^n_C\cF_{C}/\wp^{n+1}_C\cF_{C} \oplus \prod_{x \in C}
\wp^n_x\cF_{x}/\wp^{n+1}_x\cF_{x} \ra
\prod_{x \in C} \wp^n_{x,C}\cF_{x,C}/\wp^{n+1}_{x,C}\cF_{x,C}
$$
is quasi-isomorphic to the complex
$$
(A_C({\cal F}) \cap \wp^n_C\cF_{C})/(A_C({\cal F}) \cap \wp^{n+1}_C\cF_{C})
\oplus \wp^n_P\cF_{P}/\wp^{n+1}_P\cF_{P} \ra
\wp^n_{P,C}\cF_{P,C}/\wp^{n+1}_{P,C}\cF_{P,C}.
$$
\end{lemma}
This lemma is an extension of the proposition 1 above. The proves of the
both lemmas are straightforward and we will skip them.

\begin{theorem}.
Let $X$ be a projective irreducible algebraic surface over a field
     $k$, $C \subset X$ be an irreducible projective curve, and
     $P \in C$ be a smooth point on both $C$ and $X$. Let ${\cal F}$
be a torsioh free coherent sheaf on $X$.

Assume that the the surface $X - C$ is affine.
Then the following complexes are quasi-isomorphic:
\begin{quotation}

i)  the adelic complex
$$
\cF_{\eta}  \oplus  \pr_{D} \cF_D   \oplus  \pr_{x} \cF_x
 \ra
\prod_{D} (\cF_D \otimes K_D) \oplus  \prod_{x} (\cF_x \otimes K_x) \oplus
\prod_{x \in D} (\cF_x \otimes \cO_{x,D})
 \ra
$$
$$
\ra \prod_{x \in D} (\cF_x \otimes K_{x,D})
$$
for the sheaf ${\cal F}$ and

ii) the complex
$$
A({\cal F}) \oplus A_C({\cal F}) \oplus \cF_P \ra B_C({\cal F}) \oplus
B_P({\cal F}) \oplus (\cF_P \otimes
\cO_{P,C}) \ra \cF_P \otimes  K_{P,C}
$$

\end{quotation}
\end{theorem}

{\sc  Proof} will be divided into several steps. We will subsequently
     transform the adelic complex checking that every time we get a
     quasi-isomorphic complex.

     {\bf Step I}.  Consider the diagram
{\footnotesize
$$
\begin{array}{cccccccccc}
&& \pr_{D \ne C} \cF_D & \oplus & \pr_{x \in U} \cF_x &
\ra & \pr_{D \ne C} \cF_D  & \oplus & \pr_{x \in U} \cF_x & \oplus \\

&& \da && \da && \da && \da & \\

\cF_{\eta} & \oplus & \pr_{D} \cF_D  & \oplus & \pr_{x} \cF_x
& \ra &
\prod_{D} K\cF_{D} & \oplus & \prod_{x} K\cF_{x} & \oplus \\

\parallel && \da && \da && \da && \da & \\

\cF_{\eta} & \oplus &  \cF_C  & \oplus & \pr_{x \in C} \cF_x
& \ra &
(\prod_{D \ne C} K\cF_{D}/\cF_D \oplus K\cF_{C}) & \oplus &
(\prod_{x \in U} K\cF_{x}/\cF_x \oplus \prod_{x \in C} K\cF_{x}) &
\oplus
\end{array}
$$
$$
\begin{array}{cccccc}
\oplus & \prod_{x \in D \ne C} \cF_{x,D} &
\ra & \prod_{x \in D \ne C} \cF_{x,D} && \\
& \da & &  \da &&\\
\oplus & \prod_{x \in D} \cF_{x,D} &
\ra & \prod_{x \in D} K\cF_{x,D} && \\
 & \da && \da &&\\

\oplus & \prod_{x \in C} \cF_{x,C} &
\ra & \prod_{x \in D \ne C} K\cF_{x,D}/\cF_{x,D} & \oplus &
\prod_{x,C} K\cF_{x,C}
\end{array}
$$}
where $U = X - C$. The middle
     row is the full adelic complex and the first row is an exact
     subcomplex. The commutativity of the upper squares is obvious.

     The exactness follows from the  trivial
\begin{lemma}. Let $f_{1,2}: A_{1,2} \ra B$ be  homomorphisms of abelian
groups.
     The complex
     $$ 0 \ra A_1 \oplus A_2 \ra A_1 \oplus A_2 \oplus B \ra B \ra 0, $$
     where $ (a_1 \oplus a_2) \mapsto (a_1 \oplus -a_2 \oplus -f(a_1) +f(a_2)),
 (a_1 \oplus a_2 \oplus b) \mapsto (f(a_1) + f(a_2) + b)$,
     is exact.
\end{lemma}
The third row in the diagram is a quotient-complex by the subcomplex and
     we conclude that it is quasi-isomorphic to the adelic complex.

		 {\bf Step II}. We can make the same step with the adelic complex
for the sheaf ${\cal F}$ on the surface $U$. By assumption the surface
     $U$ is affine  and we get an {\em exact} complex
$$
\cF_{\eta}/A \ra \prod_{D \ne C} (\cF_D \otimes K_D)/\cF_D \oplus
\prod_{x \in U} (\cF_x \otimes K_x)/\cF_x
     \ra
$$
$$
  \prod_{\begin{array}{c}
               x \in U \\
              x \in D \ne C
               \end{array}}
     (\cF_x \otimes K_{x,D})/(\cF_x \otimes \cO_{x,D}),
$$
where $A = \Gamma(U, {\cal F})$.
\begin{lemma}. The complex
$$
0 \ra \prod_{x \in C} (\cF_x \otimes K_x)/B_x({\cal F})
     \ra \prod_{\begin{array}{c}
               x \in C \\
              x \in D \ne C
               \end{array}}
     (\cF_x \otimes K_{x,D})/(\cF_x \otimes \cO_{x,D}) \ra 0
$$
is exact.
\end{lemma}
{\sc Proof}. The injectivity follows directly from the definition of
     the ring $B_x$. The surjectivity is the local strong approximation
     around the point $x \in C$ (see \cite{P1}[\S 1],\cite{FP}[ch.4]).

     {\bf Step III}.  Take the sum of the  two complexes from
     step II. Then we have a map of the complex we got in the step I to
     this complex
{\footnotesize
$$
\begin{array}{cccccccccc}

\cF_{\eta} & \oplus &  \cF_C  & \oplus & \pr_{x} \cF_x
& \ra &
(\prod_{D \ne C} K\cF_{D}/\cF_D \oplus K\cF_{C}) & \oplus &
(\prod_{x \in U} K\cF_{x}/\cF_x  \oplus  \prod_{x \in C} K\cF_{x})
& \oplus \\

\da && \da && \da && \da && \da &  \\

\cF_{\eta}/A & \oplus & (0) & \oplus & (0)
& \ra &
\prod_{D \ne C} K\cF_{D}/\cF_D & \oplus &
(\prod_{x \in U} K\cF_{x}/\cF_x  \oplus  \prod_{x \in C} K\cF_{x}/B_x)
& \oplus
\end{array}
$$
$$
\begin{array}{cccccc}

\oplus & \prod_{x \in C} \cF_{x,C} &
\ra & \prod_{x \in D \ne C} K\cF_{x,D}/\cF_{x,D}  &
\oplus & \prod_{x \in C} K\cF_{x,C}  \\

& \da && \da && \da  \\

\oplus & (0) &
\ra & \prod_{\begin{array}{c}
x \in U \\
x \in D \\
D \ne C
\end{array}} K\cF_{x,D}/\cF_{x,D}
 \oplus
\prod_{\begin{array}{c}
x \in C \\
x \in D \\
D \ne C
\end{array}} K\cF_{x,D}/\cF_{x,D} & \oplus & (0)
\end{array}
$$}
For this map all the
     components which do not have arrows are mapped to zero. The diagram
     is commutative and the kernel of the map is equal to
$$
A \oplus \cF_C \oplus \pr_{x \in C} \cF_x  \ra
     K\cF_C \oplus \prod_{x \in C} B_x({\cal F}) \oplus
\prod_{x \in C} K\cF_x  \ra \prod_{x \in C} K\cF_{x,C}.
$$
We conclude that this complex is quasi-isomorphic to the adelic complex.

     {\bf Step IV}. Using the embedding $\cF_x \ra B_x({\cal F})$ and lemma 5
     we have an exact complex and it's embedding into the complex of the
     step III:
{\small
$$
\begin{array}{cccccccccccc}
&&&& \pr_{x \in C-P} \cF_x & \ra &&& \pr_{x \in C-P} B_x({\cal F}) & \oplus &
\pr_{x \in C-P} \cF_x &  \ra \\

&&&& \da &&&& \da && \da &  \\

A & \oplus & \cF_C & \oplus & \pr_{x \in C}  \cF_x & \ra &
K\cF_C  & \oplus & \prod_{x \in C} B_x({\cal F}) & \oplus &
\prod_{x \in C} \cF_{x,C}
  &   \ra
\end{array}
$$
$$
\begin{array}{cc}
\ra & \pr_{x \in C-P} B_x({\cal F})  \\
& \da \\
\ra &  \prod_{x \in C} K\cF_{x,C}
\end{array}
$$
}

As a result we get the factor-complex
$$
A \oplus \cF_C \oplus \cF_P \ra K\cF_C \oplus B_P({\cal F}) \oplus
     \prod_{x \in C-P} \cF_{x,C}/\cF_x \oplus \cF_{P,C} \ra
$$
$$
\ra  \prod_{x \in C-P} K\cF_{x,C}/B_x({\cal F}) \oplus K\cF_{P,C}.
$$

     {\bf Step V}. Now we need
\begin{lemma}. The complex
$$
0 \ra  (\cF_C \otimes K_C)/B_C({\cal F})
     \ra \prod_{x \in C-P}
     (\cF_x \otimes K_{x,C})/B_x({\cal F})  \ra 0
$$
is exact.
\end{lemma}
{\sc Proof}. The injectivity is again the definition of the $B_C$ and
     the surjectivity follows from the strong approximation on the
     curve $C$ (see proof of proposition 3) and  lemma 4 above.

     As a corollary we have an isomorphism
$$
\cF_C/A_C({\cal F}) \stackrel{\cong}{\ra} \prod_{x \in C-P} (\cF_x \otimes
\cO_{x,C})/\cF_x,
$$
where
$$ A_C({\cal F}) := B_C({\cal F}) \cap \cF_C. $$

Combining the isomorphisms from the lemma and its corollary into a
     single complex of length 2, we get the diagram
{\small
$$
\begin{array}{cccccccccccc}
A & \oplus & \cF_C & \oplus & \cF_P
& \ra &
K\cF_{C} & \oplus & B_P({\cal F}) & \oplus &
     (\prod_{x \in C-P} \cF_{x,C}/\cF_x \oplus \cF_{P,C})
& \ra \\

 \da && \da && \da && \da && \da && \da & \\

(0) & \oplus &   \cF_C/A_C  & \oplus & (0)
& \ra &
K\cF_C/B_C & \oplus &  (0) & \oplus & \prod_{x \in C-P} \cF_{x,C}/\cF_x
     & \ra
\end{array}
$$
$$
\begin{array}{cccc}
\ra  &
\prod_{x \in C-P} K\cF_{x,C}/B_x({\cal F}) & \oplus & K\cF_{P,C} \\
 & \da && \da \\
\ra  & \prod_{x \in C-P} K\cF_{x,C}/B_x({\cal F}) & \oplus & (0)
\end{array}
$$
}
The kernel of the map of the complexes is obviously equal to

$$
A({\cal F}) \oplus A_C({\cal F}) \oplus \cF_P \ra B_C({\cal F}) \oplus
B_P({\cal F}) \oplus (\cF_P \otimes
\cO_{P,C}) \ra  (\cF_P \otimes K_{P,C})
$$
and we arrive to the conclusion of the theorem.

{\sc Remark 3}. Sometimes we will call the complex from the theorem as the
{\em restricted} adelic complex.

\begin{lemma}.Let $X$ be a projective irreducible variety over a field $k$
and ${\cal O}(1)$ be a very ample sheaf on $X$.
Then
\begin{enumerate}
\item The following conditions are equivalent
     \begin{quotation}
    i) $X$ is a Cohen-Macaulay variety

    ii) for any locally free sheaf ${\cal F}$ on $X$ and $i < \mbox{dim}(X)$
     ~$H^{i}(X, {\cal F}(n)) = (0)$ for $n << 0$
     \end{quotation}
\item If $X$ is normal of dimension $> 1$ then for any locally free sheaf
${\cal F}$ on $X$~$H^{1}(X, {\cal F}(n)) = (0)$ for $n << 0$
\end{enumerate}
\end{lemma}
{\sc Proof} see in \cite{Har}[ch. III, Thm. 7.6, Cor. 7.8]. The last
statement is known as the lemma of Enriques-Severi-Zariski. For
dimension 2 every normal variety is Cohen-Macaulay \cite{Mat} and thus
the second claim
follows from the first one.

\begin{prop}. Let ${\cal F}$ be a locally free  coherent sheaf on the
projective irreducible surface $X$.

Assume that the local rings of the $X$ are Cohen-Macaulay and the curve
$C$ is a locally complete intersection. Then, inside the
field $K_{P,C}$, we have
$$B_C({\cal F}) \cap B_P({\cal F}) = A({\cal F}).$$
\end{prop}
{\sc Proof} will be done in several steps.

{\sc Step 1}. If we know the proposition for a sheaf ${\cal F}$ then it
is true for the sheaf ${\cal F}(nC)$ for any $n \in \Z$. Thus taking
a twist by  ${\cal O}(n)$ we can assume that $\mbox{deg}_C({\cal F}) < 0$.

{\sc Step 2}.
Now we show that $A_C({\cal F}) \cap \cF_P = (0)$. The filtrations from
lemma 3 gives the corresponding filtration of the group $A_C({\cal F})$.
Lemma 4 implies that
$$ \frac{(A_C({\cal F}) \cap \cF_P) \cap \wp^n\cF_P}{(A_C({\cal F}) \cap
\cF_P) \cap \wp^{n+1}\cF_P}
\cong \Gamma(C, {\cal F} \otimes \check{N}^{\otimes n}_{X/C}).
$$
Since $\mbox{deg}_C({\cal F}) < 0$,~$N_{X/C} = {\cal O}_X(C)\vert_C$ and
$\mbox{deg}_C(N_{X/C}) > 0$
we get that the last group is trivial.

{\sc Step 3}.
The next step is to prove the equality:
$$
B_C({\cal F}(-D)) \cap B_P({\cal F}(-D)) = A({\cal F}(-D)),
$$
where $D$ is an sufficiently ample divisor on $X$ distinct from the curve
$C$. By theorem 3 the cohomology of ${\cal F}_X(-D)$
can be computed from the complex
$$
A({\cal F}(-D)) \oplus A_C({\cal F}(-D)) \oplus \cF_P(-D) \ra
B_C({\cal F}(-D)) \oplus B_P({\cal F}(-D)) \oplus
\cF_{P,C}(-D)
$$
$$
\ra   K{\cal F}_{P,C}.
$$
Now take $a_{01} \in B_C({\cal F}(-D)), a_{02} \in B_P({\cal F}(-D))$ such
that $a_{01} + a_{02} = 0$.
They define an element $(a_{01} \oplus a_{02} \oplus 0)$ in the middle
component of the complex.
By our condition for $D$ and lemma 8 we have $H^1(X, {\cal F}_X(-D)) = (0)$
and thus there exist $a_0 \in A({\cal F}(-D)), a_1 \in A_C{\cal F}(-D),
a_2 \in \cF_P(-D)$ such that $a_{01} = a_0 - a_1, a_{02} = a_2 - a_0, 0 a_1 - a_2$.

By the second step $a_1 = a_2 \in (A_C({\cal F}(-D)) \cap \cF_P(-D)) \subset
A_C({\cal F}) \cap \cF_P = (0)$ and, consequently, we have
$a_{01} (=-a_{02}) \in A({\cal F}(-D))$.

{\sc Step 4}.
The last step is to take two distinct divisors $D, D'$ such that
$D \cap D' \subset C$. Since $C$ is a hyperplane section  we can choose for
$D, D'$ two hyperplane sections whose intersection belongs to $C$. Therefore
their ideals in the ring $A({\cal F})$ are relatively prime and
$$ A({\cal F}) = A({\cal F}(-D)) + A({\cal F}(-D')) \ni 1 = a + a' , a
\in A({\cal F}(-D)), a' \in A({\cal F}(-D')).$$
If now $b \in B_C({\cal F}) \cap B_P({\cal F})$, then $b = ba + ba'$,
where $ba \in B_C({\cal F}(-D)) \cap B_P({\cal F}(-D)), ba' \in
B_C({\cal F}(-D')) \cap B_P({\cal F}(-D'))$. We see that $b \in A({\cal F})$
by the previous step.

{\sc Remark 4}. The method we have used  cannot be applied if our
variety is not Cohen-Macaulay (by lemma 8 above). It would be interesting to
know how to extend
the result to the arbitrary surfaces $X$ and the sheaves ${\cal F}$ such
that ${\cal F}$ are locally free outside $C$. The last condition is really
necessary.
We also note that any normal surface is Cohen-Macaulay \cite{Mat}[\S 17].

{\sc Remark 5}. This proposition is a version for the reduced adelic complex
of the corresponding result for the full complex. Namely, $\A_{X,01} \cap
\A_{X,02} = \A_{X,0}$, see \cite[ch.IV]{FP}. This should be generalized
to arbitrary dimension $n$ in the following way.

Let $I, J \subset [0, 1, \dots, n]$ and
$$
\A_{X, I}({\cal F}) = (\prod_{\{codim \eta_0, codim \eta_1, \dots \} \in I}
K_{\eta_0,\eta_1,\dots}) \otimes {\cal F}_{\eta_0}) ~\bigcap ~\A_X({\cal F}).
$$
Then we have
$$
\A_{X, I}({\cal F}) \cap \A_{X, J}({\cal F}) = \A_{X, I \cap J}({\cal F})
$$
for a locally free ${\cal F}$ and a Cohen-Macaulay $X$.
\par\medskip
{\sc Example}. Let $X = {\bf P}_2 \supset C = {\bf P}_1 \ni P$. We introduce
homogenous coordinates $(x_0:x_1:x_2)$ such that  $C = (x_0 = 0); P (x_0 = x_1 = 0$ and
$U = X - C = Spec k[x,y]$ with $x = x_1/x_0, y = x_2/x_0$.
Then $k(C) = k(y/x), x^{-1}$ is the last parameter for any two-dimensional
local field $K_{Q,C}$ with $Q \ne P$. For local field $K_{P,C}$ we have
$$ K_{P,C} = k((u))((t)), u = xy^{-1}, t = y^{-1}.$$

Then we can easily compute all the rings from the complex of theorem 1 for
the sheaf ${\cal O}_X$.
$$
\begin{array}{ccccl}
&& B_P & = & k[[u]]((t)) \\
&& B_C & = & k[u^{-1}]((u^{-1}t)) \\
&& \cO_{P,C} & = & k((u))[[t]] \\
A & = & \Gamma(U, {\cal O}_X) & = & k[ut^{-1}, t^{-1}] \\
&& A_C & = & k[u^{-1}][[u^{-1}t]] \\
&& \cO_P & = & k[[u,t]]
\end{array}
$$
We can draw the subspaces as some subsets of the plane according to the
supports of the elements of the subspaces (on the plane with
coordinates $(i,j)$ for elements $u^it^j \in K_{P,C}$. Then the first
three subspaces $B_P, B_C, \cO_{P,C}$ will correspond to some halfplanes
and the subspaces $A, A_C, \cO_P$ to the intersections of them.

\section{Krichever correspondence for dimension 2}

We need the following well known result.
\begin{lemma}. Let $X$ be an projective variety, ${\cal F}$ be a coherent
sheaf on ��and $C$ be an ample divisor on $X$. If
$$
S = \oplus_{n \ge 0} \Gamma(X, {\cal O}_X(nC)),~F = \oplus_{n \ge 0}
\Gamma(X, {\cal F}(nC)),
$$
then
$$ X \cong Proj(S),~{\cal F} \cong Proj(F).$$
\end{lemma}
{\sc Proof}. Let $mC$ be a very ample divisor,$S = \oplus_{n \ge 0} S_n$ and
$S': = \oplus_{n \ge 0} S_{nm}$. Then by \cite{GD}[prop. 2.4.7]
$$Proj(S') \cong Proj(S).$$

The divisor $mC$ defines an embedding $i: X \ra {\bf P}$ to a projective
space such that $i^*{\cal O}_{\bf P}(1) = {\cal O}_X(mC)$. Let ${\cal J}_X \subset
{\cal O}_{\bf P}$ be an ideal defined by $X$. If

$$ I: = \oplus_{n \ge 0} \Gamma({\bf P}, {\cal J}_X(n)), $$
$$ A: = \oplus_{n \ge 0} \Gamma({\bf P}, {\cal O}_{\bf P}(n)), $$
then $I \subset A$ and by \cite{GD}[prop. 2.9.2]
	$$ Proj(A/I) \cong X.$$

	We have an exact sequence of sheaves
	$$0 \ra {\cal J}_X(n) \ra {\cal O}_{\bf P} \ra {\cal O}_X(n) \ra 0,$$
	which implies the sequence
	$$ 0 \ra \bigoplus\limits_{n \ge 0} \Gamma({\cal J}_X(n)) \ra
        \bigoplus\limits_{n \ge 0}
	\Gamma({\cal O}_{\bf P}(n)) \ra \bigoplus\limits_{n \ge 0} \Gamma({\bf P},
	{\cal O}_X(n)) \ra  \bigoplus\limits_{n \ge 0} H^1({\bf P}, {\cal J}_X(n)).$$
	Here the last term is trivial for sufficiently large $n$. The first
	three terms are equal correspondingly to $I, A$ and $S'$. It means that
	the homogenous components of $A/I$ and $S' \supset A/I$ are equal for
	sufficiently big degrees.

	By \cite{GD}[prop. 2.9.1]
	$$Proj(A/I) \cong Proj(S'),$$
	and combining everything together we get the statement of the lemma.
        The statement concerning the sheaf ${\cal F}$ can be proved along
        the same line.
\par\medskip

Now we move to the case of algebraic surfaces. The corresponding data
has the following

{\sc Definition 4}.
 $$
 \begin{array}{lll}
 {\cal M}_2 & := & \{ X, C, P, (z_1, z_2),  {\cal F}, e_P \}\\
 X & & \mbox{projective irreducible surface}~/k \\
 C \subset X & & \mbox{projective irreducible curve}~/k \\
 P \in C && \mbox{a smooth point on }~X~\mbox{and}~C \\
 z_1, z_2 && \mbox{formal local parameter at}~P~\mbox{such that} \\
        && (z_2 = 0) = C ~\mbox{near}~P \\
 {\cal F} && \mbox{torsion free rank}~$r$~\mbox{sheaf on}~X \\
 e_P && \mbox{a trivialization of}~{\cal F}~\mbox{at}~P
 \end{array}
 $$
Then we have
$$
\cO_{X,P}= k[[z_1, z_2]],~K_{P,C} = k((z_1))((z_2)),
$$
$$
\cF_P  = \cO_Pe_P = \cO_P^{\oplus r}.
$$
For the field $K = k((z_1))((z_2))$ we have the following filtrations and
subspaces:
$$
\begin{array}{lll}
K_{02} & = & k[[z_1]]((z_2)), \\
K_{12} & = & k((z_1))[[z_2]],   \\
K(n) & = & z_2^n K_{12}.
\end{array}
$$
Taking the direct sums we introduce the subspaces $V_{02}, V_{12},V(n)$ of
the space $V = K^{\oplus r}$.
\begin{theorem}. Let~$C$ be a hyperplane section on the surface $X$. Then
there exists a canonical map
$$
\Phi_2 : {\cal M}_2 \ra \{\mbox{vector subspaces}~B \subset K, W \subset V \}
$$
such that
\begin{quotation}

i) for all $n$ the complexes
 $$
 \frac{B \cap K(n)}{B \cap K(n+1)} \oplus \frac{K_{02} \cap K(n)}{K_{02} \cap
  K(n+1)} \ra \frac{K(n)}{K(n+1)}
$$
 $$
 \frac{W \cap V(n)}{W \cap V(n+1)} \oplus \frac{V_{02} \cap V(n)}{V_{02} \cap
  V(n+1)} \ra \frac{V(n)}{V(n+1)}
$$

are Fredholm of index $\chi(C, {\cal O}_C)  + nC.C$ and
$\chi(C, {\cal F}\vert_C)  + nC.C$, respectively

ii) the cohomology of  complexes
$$
(B \cap K_{02}) \oplus (B \cap K_{12}) \oplus (K_{02}) \cap K_{12})  \ra
B \oplus K_{02} \oplus K_{12} \ra K
$$
$$
(W \cap V_{02}) \oplus (W \cap V_{12}) \oplus (V_{02}) \cap V_{12})  \ra
W \oplus V_{02} \oplus V_{12} \ra V
$$
are isomorphic to $H^{\cdot}(X, {\cal O}_X)$ and $H^{\cdot}(X, {\cal F})$,
respectively

iii)if $(B, W) \in \mbox{Im}~\Phi_2$ then $B \cdot B \subset B,
B \cdot W \subset W$

iv) for all $n$ the map
$$
(C, P, z_1\vert_C, {\cal F}(nC)\vert_C, e_P(n)\vert_C) \mapsto
$$
$$
\mapsto (\frac{B \cap K(n)}{B \cap K(n+1)}  \subset  \frac{K(n)}{K(n+1)} k((z_1)),
$$
$$
\frac{W \cap V(n)}{W \cap V(n+1)} \subset  \frac{V(n)}{V(n+1)} k((z_1))^{\oplus r}~)
$$
coincides with the map $\Phi_1$.

v) let the sheaf ${\cal F}$ be locally free and the surface $X$ be Cohen-
Macaulay. If $m,m^{\prime}
\in {\cal M}_1$ and $\Phi_2(m) = \Phi_2(m^{\prime})$ then $m$ is isomorphic
to $m^{\prime}$
\end{quotation}
\end{theorem}

{\sc Proof}. If $ m = (X, C, P, (z_1, z_2),  {\cal F}, e_P) \in {\cal M}_2$ then
to define the map $\Phi_2$ we put
$$
\begin{array}{ccl}
B &  = & B_C({\cal O}_X), \\
W &  = & B_C({\cal F}), \\
\Phi_2(m) & = & (B, W).
\end{array}
$$
Since we have the local coordinates $z_{1,2}$ and the trivialization $e_P$
the subspaces $B$ and $W$ will belong to the space $k((z_1))((z_2))$ exactly
as in the case of dimension 1 considered above.

We note that our condition on the curve $C$ implies that $C$ ia  a Cartier
divisor and the surface $X - C$ is affine.

The property i) follows
from lemma 4, the property ii) follows from theorem 3 and the general
theorem of adelic theory: the cohomologies of the adelic complex of a sheaf
${\cal F}$ are equal to the cohomologies $H^{\cdot}(X, {\cal F})$  of the
sheaf ${\cal F}$~\cite{P1, Be, H, FP}. The property
iii) is trivial again, to get iv) one needs again to apply lemma 4 and
to get v) it is enough to use proposition 4 and lemma 9. They show that given
a point $(B, W) \in {\cal M}_2$ such that $(B, W) = \Phi_2(m)$ we can
reconstruct the data $m$ up to an isomorphism.

{\sc Remark 6}. The property v) of the theorem cannot be extended to the
arbitrary torsion free sheaves on $X$. We certainly cannot reconstruct
such sheaf if it is not locally free outside $C$. Indeed, if ${\cal F},
{\cal F}^{\prime}$ are two sheaves and there is a monomorphism
${\cal F}^{\prime} \ra {\cal F}$ such that ${\cal F}/{\cal F}^{\prime}$ has
support in $X - C$ then the restricted adelic complexes for the
sheaves ${\cal F}, {\cal F}^{\prime}$ are isomorphic.

{\sc Remark 7}. A definition of the map $\Phi_n$ for all $n$ was
suggested in \cite{O}. It has the properties that correspond to the
properties i) - v) of the theorem. Also the proofs in \cite{O} has  the
advantage: they are direct and don't use the general adelic machinery.

\begin{center}
{\bf Appendix 1}
\end{center}
Here we show how to deduce from the Lax form of the KP hierarchy for
pseudo-differential operators $L$ the classical KP and KdV equations.

Let $L = \partial +  u_1\partial^{-1} + u_2\partial^{-2} + \dots$,
where $u_m = u_m(x, t_1, t_2, \dots) \in k[[x, t_1, t_2, \dots]]$ for $m \ge 1$.
If we denote $\partial/\partial t_n$ as $\partial_n$ then the KP hierarchy
has the following form
$$
\partial_n L = [(L^n)_{+}, L].
$$
This gives us for any $n$ an infinite sequence of differential equations for
the functions $u_m$. Denote by $u', u'', \dots$ the subsequent derivatives
by $x$.

First for $n = 1$, we get
$$
\partial_1 u_m =u'_m~\mbox{for all}~m \ge 1
$$
This means that we can take $t_1 = x$.

Now we write down the first two equations for $n = 2$ and the first
equation for $n = 3$.
$$
\partial_2 u_1 = u''_1 + 2u'_2\eqno (1)
$$
$$
\partial_2 u_2 = u''_2 + 2u'_3 + 2u_1u'_1\eqno (2)
$$
$$
\partial_3 u_1 = u'''_1 + 3u''_2 +  3u'_3 + 6u_1u'_1\eqno (3)
$$
Let us introduce the new notations: $u = u_1(x, y, t)$ with $y = t_2,
t = t_3$. Also we use the standard notations $u_t, u_y, u_{yy},\dots$ for
derivatives.

We can eliminate $u'_3$ from equations (2) and (3) and then we get
$$
2u_t - 2u''' - 6uu' = 3(u''_2 + u_{2y}).\eqno (4)
$$
Differentiating this equation by $x$ we have
$$
(2u_t - 2u''' - 6uu')' = 3(u'''_2 + u'_{2y})\eqno (5).
$$
>From (1) we  find
$$
u'''_2 = 1/2(u''_y - u''''),~u'_{2y} =  1/2(u_{yy}  - u''_y)
$$
and inserting these expressions into (5) we finally get  the KP equation
$$
(4u_t - u''' - 12uu')' = 3u_{yy}.
$$
In the space $E'$ there is an invariant submanifold defined by condition
$(L^2)_{-} = 0$ \cite{P4}[Lemma 2]. Let us look what this condition means
for the operator $L$ as above. We have
$$
L^2 = \partial^2 + (2u_2 + u'_1)\partial^{-1} + \dots
$$
and thus in our notations
$$
2u_2 + u' = 0.\eqno (6)
$$
Combining (1) and (6) we see that $u_{2y} = 0$. Together with (4) this
implies $3u''_2 = 2u_t - 2u''' - 6uu'$. Using once more (6) we get at last
the KdV equation
$$
4u_t - 7u''' - 12uu' = 0.
$$
The numerical coefficients here are not essential since one can get all
possible values by rescaling of $x, t$ and $u$.

These computations and the choice of variables look quite artificial
comparing with the previous constructions. It is an interesting independent
problem how to deduce the KP and KdV equations from the Lax operator
equations in a more conceptual way.

Now,  we show how to do that for the classical KdV equation.
Comparing with the previous deduction this one can be done in a purely
formal way, without any "tricks".

Again, we order that our equation will satisfy the constraint condition
$(L^2)_- = 0$. For the operator $L$ as
above we have
$$
L^2 = \partial^2 + (2u_2 + u'_1)\partial^{-1} + (2u_3 + u_1^2 + u'_2)\partial^2 + \dots
$$
and thus in our notations
$$
2u_2 + u' =0,
2u_3 + u_1^2 + u'_2  = 0
$$
and so on. This means that we can compute all $u_m$ for $m> 1$ starting from $u_1$ and it's
derivatives. Under this condition the first non-trivial equation of the
KP hierarchy is the following one
$$
\partial_3L = [(L^3)_+, L].
$$
and taking the coefficient nearby $\partial^{-1}$ we have
$$
u'''_1 +3u''_2 +3u'_3 +6u_1u'_1 = 0.
$$

Substituting the given above expressions for $u_2$ and $u_3$ we get at last the
KdV equation
$$
4u_t - 7u''' - 12uu' = 0.
$$

{\bf Problem.} Let us go to the case of  dimension two. Then we have  the constraint conditions of the type
$$
(L^mM^n)_- = 0
$$
and the particular components of  hierarchy.
Can we deduce some concrete equations using the same way as above ? Certainly, we have more possibilities. The
first question is how many initial functions $u_{mn}$ one has to use so to generate all other coefficients ? For dimension 1 only $u_1$ was enough !

{\bf Remark.} For higher dimensions it is crucial to transform the whole system to
avoid the asymmetry of the variables. It means that we have to consider not
the ring $P$ alone but at least the direct sum $P\oplus P'$ where $P'$ has interchanged variables, namely $y, x$ instead of $x, y$. As a goal one can hope to get in this way the equations of the plane hydrodynamics.
It is known they have infinitely many conservation laws.

\begin{center}
{\bf Appendix 2}
\end{center}
Here, we present a well-known construction of semi-infinite monomes in an infinite-dimensional vector space.
Let $V$ be a vector space over a fied $k$ and $V_n, n \in \Z$ be an exhausted increasing filtration in $V$. The example
one has to have in mind is when $V = k((z)), V_n = z^{-n}k[[z]]$. If $V$ would be a finite-dimensional space,the
$k$-dimensional vector subspaces $W \subset V$ can be described by the elements $ \bigwedge^k(W) \subset \bigwedge^k(V)$.
In the infinite-dimensional case we assume that
$$
\mbox{dim} W \cap V_n  = k + n~ \mbox{for large}~ n.
$$
In the case of $V = k((z))$ this $k$ is exactly the index of subspace $W$ from the Sato Grassmanian.
In general case for such $n$, we have the diagram
$$
\begin{diagram}
\node{W\cap V_n}         \arrow{e,t}{} \arrow{s,l}{}
\node{W\cap V_{n + 1}}         \arrow{s,r}{}\\
\node{ V_n}
\arrow{e,t}{} \node{ V_{n + 1}}
\end{diagram}
$$
We denote the one-dimensional space $(W\cap V_{n + 1})/(W\cap V_n)$ as $1_n$. The diagram induces
isomorphism $1_n \rightarrow V_{n + 1}/V_n$.

For any exact sequence $0 \rightarrow V \rightarrow V' \rightarrow V'/V \rightarrow 0$ with
one-dimensional space $V'/V$, there is a canonical map $\bigwedge^{n+1}(V') \rightarrow \bigwedge^n(V)\otimes (V'/V)$.
If $v = v_1 \wedge v_2 \wedge \dots v_{n+1} \in \bigwedge^{n+1}(V')$ then the image of $v$ is zero if either
all $v_i$ belong to $V$ or at least two $v_i$ belong to $V' - V$. If exactly one $v_i$ belong to $V'- V$ then
the image is $(-1)^{n-i}v_1 \dots \wedge v_{i-1}  \wedge v_{i+1}\wedge \dots v_{n+1}\otimes v_i \mbox{mod}V$.

We get the new diagram
$$
\begin{diagram}
\node{\bigwedge^{k+n+1}(W\cap V_{n + 1})} \arrow{e,t}{} \arrow{s,l}{}
\node{\bigwedge^{k+n}(W\cap V_n)\otimes 1_n}         \arrow{s,r}{}\\
\node{\bigwedge^{k+n+1}(V_{n + 1})}
\arrow{e,t}{} \node{\bigwedge^{k+n}(V_n)\otimes V_{n+1}/V_n}
\end{diagram}
$$
For any two subspaces $A, B \subset V$, we say they are commesurable  iff the intersection
$A \cap B$ has finite codimension in both of them. In this case, we define
$$
(A\mid B) = \mbox{det}(A/A\cap B)^{-1}\otimes \mbox{det}(B/A\cap B),
$$
where $\mbox{det}(A) = \bigwedge^{\mbox{dim}(A)}(A)$ and $(A)^{-1}$ is the space dual to $A$.
Note that all the spaces $(A\mid B)$ are one-dimensional. There are the canonical
isomorphisms $(A\mid B) \cong (B\mid A)^{-1}$ and $(A\mid C) \cong (A\mid B)\otimes (B\mid C)$ which
we will treat as equalities.
These rules allows us to write
$$
V_{n+1}/V_n = (V_n\mid V_{n+1} ) = (V_n\mid V_0 )( V_{n+1}\mid V_0)^{-1},
$$
$$
1_n  = (V_n\mid V_{n+1} ) = (V_n\mid V_0 )( V_{n+1}\mid V_0)^{-1}~\mbox{for large}~n.
$$
Thus, we can rewrite the diagram as
$$
\begin{diagram}
\node{\bigwedge^{k+n+1}(W\cap V_{n + 1})( V_{n+1}\mid V_0))} \arrow{e,t}{} \arrow{s,l}{}
\node{\bigwedge^{k+n}(W\cap V_n)\otimes (V_n\mid V_0 )}         \arrow{s,r}{}\\
\node{\bigwedge^{k+n+1}(V_{n + 1})( V_{n+1}\mid V_0))}
\arrow{e,t}{} \node{\bigwedge^{k+n}(V_n)\otimes (V_n\mid V_0 )}
\end{diagram}
$$
At last, the space of semi-infinite monomes (of index $k$) can be defined as projective
limit respect these maps:
$$
\bigwedge^{k+\frac{\infty}{2}}(V) = \mbox{lim}_n\bigwedge^{k+n}(V_n)\otimes (V_n\mid V_0 )
$$
and we can attach to the subspace $W$ the line
$$
[W] =  \mbox{lim}_n\bigwedge^{k+n}(W\cap V_n)\otimes (V_n\mid V_0 ).
$$

Using these constructions we can present the Sato Grassmanian as a disjoint union of connected components
$$
Gr(V) = \coprod_k Gr_k(V)
$$
and every component has a projective embedding
$$
Gr_k(V) \rightarrow  {\bf P}(\bigwedge^{k+\frac{\infty}{2}}(V))
$$
where $W \mapsto [W]$. Note that the vector space $\wedge^{k+\frac{\infty}{2}}(V)$ does depend on the choice of the subspace $V_0$ from the filtration but it's projectivization doesn't.

\par\medskip
\begin{flushleft}
Steklov Mathematical Institute\\
Gubkina str. 8\\
117966 Moscow GSP-1\\
Russia\\
e-mail: parshin@mi.ras.ru
\end{flushleft}

\end{document}